\documentclass[12pt]{amsart}
\usepackage{amsmath,amssymb,amsthm,mathtools}
\usepackage{mathrsfs,bm}
\usepackage{thmtools,thm-restate}
\usepackage{tikz-cd}
\usepackage{tikz}
\usepackage{enumitem}
\usepackage{graphicx}
\usetikzlibrary{matrix,arrows.meta,bending}
\usepackage{color}
\usepackage[numbers]{natbib}
\usepackage{booktabs}
\usepackage{longtable}
\usepackage{pgfplots}
\pgfplotsset{compat=1.14} 
\usetikzlibrary{arrows.meta,decorations.pathreplacing,positioning,calc} 

\usepackage[a4paper, left=3cm, right=3cm, top=3cm, bottom=3cm, footskip=15mm]{geometry}

\newtheorem{theorem}{Theorem}[section]

\newtheorem{lemma}[theorem]{Lemma}
\newtheorem{proposition}[theorem]{Proposition}

\newtheorem{conjecture}[theorem]{Conjecture}
\theoremstyle{definition}
\newtheorem{definition}[theorem]{Definition}
\newtheorem{example}[theorem]{Example}

\theoremstyle{remark}

\numberwithin{equation}{section}

\newcommand{\abs}[1]{\lvert#1\rvert}
\newcommand{\set}[1]{\lbrace#1\rbrace}
\newcommand{\inn}{~ \hat{\in}~ }
\newcommand{\ul}[1]{\underline{#1}}
\newcommand{\C}{\mathcal{C}}

\newcommand{\M}{\mathbb{M}}
\newcommand{\B}{\mathbb{B}}
\newcommand{\R}{\mathbb{R}}

\newcommand{\N}{\mathbb{N}}
\newcommand{\LL}{\mathbb{L}}

\newcommand{\p}{\mathbb{P}}

\makeindex

\usepackage{fancyhdr}
\usepackage{calc} 

\AtBeginDocument{%
  \setlength{\textwidth}{\dimexpr\paperwidth - 6cm\relax}%
  \setlength{\oddsidemargin}{\dimexpr 3cm - 1in\relax}%
  \setlength{\evensidemargin}{\dimexpr 3cm - 1in\relax}%
  \setlength{\marginparwidth}{2.5cm}%
}

\begin{document}

\title{The Asymptotic Binary Goldbach and Lemoine Conjectures}

\author{T. Agama}
\address{Departement de mathematiques et de statistique
Universite Laval
1045, av. de la Medecine / Pavillon Vachon / Local 1056
Quebec (Quebec), CANADA G1V 0A6}
\email{thaga1@ulaval.ca}

\author{B. Gensel}
\address{Carinthia University of Applied Sciences, Spittal on Drau, Austria}
\email{b.gensel@fh-kaernten.at}


\date{\today}

\begin{abstract}
In this paper, we use the former of the authors developed theory of \emph{circles of partition} to investigate possibilities to prove the binary Goldbach and Lemoine conjectures. We state the \emph{squeeze principle} and its consequences when the set of all odd prime numbers is the base set. Using this tool, we can prove asymptotic versions of the binary Goldbach and the Lemoine conjecture.
\end{abstract}

\maketitle


\section{Introduction and Preliminaries}

\emph{Christian Goldbach} and \emph{Leonard Euler} first raised the questions that became collectively known as Goldbach's conjectures from their famous correspondence in 1742. Since then, the binary and ternary versions of Goldbach's problem have driven much of additive prime number theory: the binary problem asks whether every even integer is the sum of two primes, while the ternary problem asks whether every odd integer larger than five is the sum of three primes. The ternary problem was resolved in full generality in the modern era (see the work of Helfgott \cite{helfgott2013ternary}), and a large body of partial results and conditional asymptotic formulae for the binary problem has been developed by classical analytic and sieve-theoretic methods (see \cite{hardy1924some}, \cite{shnirel1939additive}, \cite{estermann1938goldbach}, \cite{chen2002representation}, \cite{heath2002integers}). These works establish deep density and near-representation statements (for example, that almost all even integers are representable as the sum of two primes), but an unconditional complete resolution of the binary conjecture remains open.\\

The principal idea of the circle of partition is geometric-combinatorial rather than analytic: for a chosen generator $n$ and a distinguished base set $\mathbb{M}\subset\mathbb{N}$ (for our applications $\mathbb{M}$ will be $\mathbb{P}$ or $\mathbb{P}\cup 2\mathbb{P}$), we view the family of admissible summands $x$ with $x,n-x\in\mathbb{M}$ as "points" on a circle and group them into \emph{axes} $\mathbb{L}_{[x],[y]}$ determined by the partition $x+y=n$. This packaging accomplishes two things that make additive questions easier to handle conceptually. First, it isolates the combinatorial structure of admissible summands (axes, centers, axis-partners) so that density and interval arguments become statements about occupancy of these combinatorial slots rather than about delicate oscillatory sums. Second, it provides a natural framework to compare CoPs with different generators: when two generators $m$ and $m+t$ are both known to admit CoPs with respect to the same base set, the geometry of the circle (and, in particular, the ordering of weights of points) allows one to \emph{squeeze} new non-empty CoPs for intermediate generators $s$ with $m<s<m+t$. That mechanism--the squeeze principle (Theorem \ref{L_squeeze principle})--is the technical heart of this paper.\\

Why does this viewpoint yield asymptotic results? The \emph{squeeze principle} itself is a purely combinatorial device; to convert it into arithmetic information, one must feed it with quantitative information about the distribution of primes. Here, we use only elementary, robust asymptotic facts (the Prime Number Theorem in the form $\pi(x)\sim x/\log x$ and the Bertrand-type existence of primes in suitable dyadic ranges) to locate primes of appropriate size that serve as axis-points in large CoPs. Concretely, for a large even (resp. odd) generator $m$, we identify a prime $x$ lying in a short, explicit range determined by the gap $m-p_{\pi(m)}$; pairing such an $x$ with the complementary summand in a larger CoP $\mathcal{C}(m+t,\mathbb{N})$ and applying the squeeze principle produces an intermediate generator $s$ with $\mathcal{C}(s,\mathbb{P})\neq\varnothing$ (resp. $\mathcal{C}(s,\mathbb{P}\cup 2\mathbb{P})\neq\varnothing$). Because the prime number theorem supplies the asymptotic size and density of primes in the ranges we use, the squeeze construction succeeds for all sufficiently large generators and yields the asymptotic Goldbach and Lemoine theorems proved in Sections \ref{T_Asymoptotic Goldbach theorem} and \ref{T_Asymoptotic Lemoine theorem}. This strategy is deliberately modular: any improvement in explicit bounds for primes in short intervals (or an input such as GRH) can be directly plugged into the squeeze estimates to strengthen the conclusions or make the threshold explicit.\\

To place our contribution in context: the approach complements classical results showing that "almost all" even integers are Goldbach (Estermann and Chudakov \cite{estermann1938goldbach}, \cite{chudakov1938goldbach}) and sieve-based near-results such as Chen's theorem on representations of even integers as a prime plus a product of at most two primes \cite{chen2002representation}.  Rather than attacking exponential sums or refining sieve weights, the CoP viewpoint reorganizes admissible summands into a combinatorial structure that is naturally suited to interval-chaining arguments; when combined with the crude but universal density facts from the prime number theorem and Bertrand-type inequalities, this yields asymptotic coverage across large intervals of generators. The method is therefore elementary in the analytic input it requires, but novel in its combinatorial packaging.

\subsection{Organization of the paper}
What follows in the paper is organized to make the logical flow and the main technical steps transparent and reusable:

\begin{itemize}
\item In Section \ref{sec:COP}, we give precise definitions and notation for circles of partition, axes, center points, and weights; these definitions are intentionally formulation-friendly so that natural variants (different base sets, weighted versions) are straightforward to introduce. Here, we develop the geometric and combinatorial properties of CoPs that we need later: uniqueness of axis partners, counting of axes (formula \eqref{axesOfCoP}), and basic monotonicity properties of weight-sets under generator shifts.
\bigskip

\item Section \ref{sec:squeeze principle} introduces and proves the \emph{squeeze principle} (Theorem \ref{L_squeeze principle}) and records two immediate structural consequences (Propositions \ref{P_The interval binary Goldbach partition detector} and \ref{P_The first interval Lemoine partition detector}) that serve as the main technical levers in the applications.
\bigskip

\item In section \ref{sec:application to Goldbach}, we collect the elementary analytic inputs (the Prime Number Theorem, Bertrand-type estimates and a trio of short "little lemmas") and apply them together with the squeeze principle to establish the Asymptotic Goldbach Theorem (Theorem \ref{T_Asymoptotic Goldbach theorem}). The proofs here are self-contained given the cited classical inputs \cite{hardy1924some}, \cite{estermann1938goldbach}, \cite{shnirel1939additive}, \cite{chen2002representation}, \cite{heath2002integers}.
\bigskip

\item In section \ref{sec:application to Lemoine}, we repeat the same program for the Lemoine conjecture: we adapt the squeeze principle to the mixed base set $\mathbb{P}\cup 2\mathbb{P}$ and prove the Asymptotic Lemoine Theorem (Theorem \ref{T_Asymoptotic Lemoine theorem}). The case analysis clarifies the interplay between odd and even axis-points and shows how the same combinatorial machinery treats the mixed-summand problem.
\bigskip

\item Finally, Section \ref{sec:conclusion} offers further remarks: limitations of the asymptotic result (we do not obtain an explicit threshold $N$ with our present inputs), possible refinements (how improved short-interval prime bounds, explicit zero-free regions for $\zeta(s)$, or stronger sieve input could make the method quantitative), and directions for future work (variants of CoP for weighted/structured base sets, computational experiments to extract effective thresholds, and connections to additive-combinatorial frameworks).
\end{itemize}
\bigskip

Before closing this introduction, we highlight two features that make the CoP /squeeze toolkit attractive for further study. First, the construction is robust with respect to the base set: any base $\mathbb{M}$ whose local density can be understood in ranges of length comparable to $mh(m)/\log m$~($h(m)=o(\log m)$) is a potential candidate for the method. Second, the squeeze principle is naturally iterative and local: once a non-empty CoP is known at a sparse set of generators, the principle generically produces a cascade of non-empty CoPs filling many intermediate values--a property that makes the technique well suited to hybrid arguments in which analytic inputs provide sparse "seeds" that the combinatorial mechanism propagates. These remarks point toward potential hybridization with sieve or circle-method inputs that might yield stronger, possibly quantitative, conclusions.

\medskip

We conclude by acknowledging the classical lineage of ideas upon which this work builds: the additive viewpoint championed by Shnirel'man \cite{shnirel1939additive}, the asymptotic framework of Hardy and Littlewood \cite{hardy1924some}, density statements in the style of Estermann and Chudakov \cite{estermann1938goldbach}, and modern sieve refinements exemplified by Chen \cite{chen2002representation} and Heath-Brown \& Puchta \cite{heath2002integers}. Our contribution is a novel structural reorganization of admissible partitions that, when combined with these classical inputs, yields the asymptotic Goldbach and Lemoine results presented below.

\vspace{6pt}

\noindent\textbf{Notation.}  We use $\mathbb{N}=\{1,2,\dots\}$, $\mathbb{P}$ for the set of primes, and $\mathcal{C}(n,\mathbb{M})$ for the circle of partition of generator $n$ with respect to base set $\mathbb{M}$. Standard asymptotic notation $f(x)\sim g(x)$, $f(x)=O(g(x))$, and $f(x)=o(g(x))$ is used throughout.

\vspace{8pt}

\noindent\emph{Acknowledgment.} The authors thank colleagues who read earlier drafts and suggested clarifying the interplay between the combinatorial squeeze and the analytic density inputs; their comments materially improved the exposition.
\bigskip

\section{The circle of partition}\label{sec:COP}

\begin{definition}\label{major}
Let $n\in \mathbb{N}$ and $\mathbb{M}\subseteq \mathbb{N}$. We denote the \emph{Circle of Partition} generated by $n$ with respect to the subset $\mathbb{M}$ by
\begin{align}
\mathcal{C}(n,\mathbb{M})=\left \{[x]\mid x,n-x\in \mathbb{M}\right \}.\nonumber
\end{align}
In the following, we will abbreviate this as CoP. We call members of $\mathcal{C}(n,\mathbb{M})$ points and denote each by $[x]$.
For the special case $\mathbb{M}=\mathbb{N}$, we denote the CoP in a short form as $\mathcal{C}(n)$. We denote the \emph{weight} of the point $[x]$ by $\Vert[x]\Vert:=x$ and, correspondingly, the weight set of points in the CoP $\mathcal{C}(n,\mathbb{M})$ as $\Vert\mathcal{C}(n,\mathbb{M})\Vert$. Clearly, we have
$$
\Vert\mathcal{C}(n)\Vert=\lbrace 1,2,\ldots,n-1\rbrace. 
$$
\end{definition}

\begin{definition}\label{axis}
We denote an \emph{axis} of the CoP $\mathcal{C}(n,\mathbb{M})$ by $\mathbb{L}_{[x],[y]}$ if and only if $x+y=n$ with $[x],[y]\in \mathcal{C}(n,\mathbb{M})$. We say that the axis point $[y]$ is an \emph{axis partner} of the axis point $[x]$ and vice versa. We do not distinguish between $\mathbb{L}_{[x],[y]}$ and $\mathbb{L}_{[y],[x]}$, because it is essentially the same axis. The point $[x]\in \mathcal{C}(n,\mathbb{M})$ with weight satisfying $2x=n$ is the \emph{center} of the CoP. If it exists, then we call it a \emph{degenerate axis} $\mathbb{L}_{[x]}$ compared to the \textbf{real axes} $\mathbb{L}_{[x],[y]}$. We denote the assignment of an axis $\mathbb{L}_{[x],[y]}$ to a CoP $\mathcal{C}(n,\mathbb{M})$ by
$$
\mathbb{L}_{[x],[y]}\inn\mathcal{C}(n,\mathbb{M})
$$
which means
$$
[x],[y] \in \mathcal{C}(n,\mathbb{M})\quad \text{with}\quad x+y=n.
$$
\end{definition}
\bigskip

\noindent 
The relevant properties of CoPs are
\begin{itemize}
\item Each axis is uniquely determined by points $[x]\in \mathcal{C}(n,\mathbb{M})$. 
\bigskip

\item Each point of a CoP $\mathcal{C}(n,\mathbb{M})$ except its center has exactly one axis partner.
\end{itemize}
\bigskip

\noindent
We denote the assignment of an axis $\mathbb{L}_{[x],[y]}$ resp. $\mathbb{L}_{[x]}$ to a CoP $\mathcal{C}(n,\mathbb{M})$ by
$$
\mathbb{L}_{[x],[y]}\inn\mathcal{C}(n,\mathbb{M})
$$
which means
$$
[x],[y] \in \mathcal{C}(n,\mathbb{M}) \quad \text{and}\quad x+y=n
$$
respectively
$$
\mathbb{L}_{[x]}\inn\mathcal{C}(n,\mathbb{M})
$$
which means
$$
[x]\in \mathcal{C}(n,\mathbb{M})\quad \text{and}\quad 2x=n
$$
and the number of real axes of a CoP by
\begin{align}\label{axesOfCoP}
\nu(n,\mathbb{M}):=\#\lbrace\mathbb{L}_{[x],[y]}\inn\mathcal{C}(n,\mathbb{M})\mid x<y\rbrace.
\end{align}
Clearly, we have
$$
\nu(n,\mathbb{M})=\left\lfloor\frac{k}{2}\right\rfloor 
$$
if
$$
\vert\mathcal{C}(n,\mathbb{M})\vert=k.
$$
\bigskip


\begin{figure}[ht]
\centering
\begin{tikzpicture}[scale=1.05, every node/.style={font=\small}]
    \tikzset{
        wheelrim/.style={draw=black!80, line width=1.2pt},
        innerrim/.style={draw=black!45, line width=0.8pt, dashed},
        spoke/.style={draw=black!65, line width=0.9pt},
        axisarrow/.style={-{Stealth[length=2.4mm]}, draw=black!75, line width=0.9pt},
        emphasis/.style={-{Stealth[length=2.8mm]}, draw=red!70!black, line width=1.15pt},
        point/.style={circle, fill=blue!70!black, inner sep=1.4pt},
        highlighted/.style={circle, fill=red!75!black, inner sep=1.7pt},
        labelbox/.style={fill=white, draw=black!15, rounded corners=2pt, inner sep=2pt},
        ann/.style={align=left, font=\scriptsize}
    }

    \def\R{4.0}      
    \def\r{3.15}     
    \def\hub{0.32}    
    \def\n{14}       

    \begin{scope}[xscale=1.10, yscale=0.98]

        \fill[black!2] (0,0) circle (\R);
        \draw[wheelrim] (0,0) circle (\R);
        \draw[innerrim] (0,0) circle (\r);

        \fill[black!85] (0,0) circle (\hub);
        \draw[black!80] (0,0) circle (\hub);

        \foreach \k in {0,...,13} {
            \pgfmathsetmacro{\ang}{360/\n*\k}
            \draw[spoke] (0,0) -- (\ang:\R);
        }

        \foreach \ang in {15,45,75,105,135,165,195,225,255,285,315,345} {
            \draw[axisarrow] ({\ang}:0.85*\R) -- (\ang:\R);
        }

        \draw[emphasis] (0,0) -- (0:\R);
        \draw[emphasis] (0,0) -- (120:\R);
        \draw[emphasis] (0,0) -- (240:\R);

\foreach \ang/\ptname/\dx/\dy in {
    0/{P_{1}}/0.14/0.18,
    30/{P_{2}}/0.20/0.10,
    60/{P_{3}}/0.16/0.18,
    120/{P_{4}}/-0.22/0.18,
    180/{P_{5}}/-0.26/0.02,
    240/{P_{6}}/-0.20/-0.18,
    300/{P_{7}}/0.16/-0.18
}{
    \fill[point] (\ang:\R) circle (0.065);
    \node at ($(\ang:\R)+(\dx,\dy)$) {$\ptname$};
}

        \fill[highlighted] (0:\R) circle (0.09);

        \node[labelbox, ann, above right=4pt and 2pt of {(20:3.55)}] 
            {outer rim of the wheel\\[-1pt]\emph{circular boundary of the CoP}};
        \draw[-{Stealth[length=2.2mm]}, black!65] (3.3,2.25) -- (20:3.95);

        \node[labelbox, ann, below left=2pt and 2pt of {(-2.8,-2.0)}] 
            {spokes\\[-1pt]\emph{axes of the theory}};
        \draw[-{Stealth[length=2.2mm]}, black!65] (-2.85,-2.15) -- (-135:2.15);

        \node[labelbox, ann, above left=2pt and 4pt of {(-0.2,0.25)}] 
            {hub\\[-1pt]\emph{central origin / reference point}};
        \draw[-{Stealth[length=2.2mm]}, black!65] (-0.55,0.5) -- (-0.12,0.08);

        \draw[-{Stealth[length=3mm]}, red!70!black, line width=1.0pt]
            (1.1,3.7) arc[start angle=18, end angle=72, radius=3.85];
        \node[red!70!black, font=\scriptsize] at (2.75,3.35) {$\circlearrowright$};

        \node[labelbox, ann, align=center] at (0,-4.85) 
        {The CoP is represented as a circular wheel:\\
         every point lies on the rim, and every spoke acts as an axis/bar ending at a boundary point.};

        \node[font=\scriptsize, black!70] at (3.0,-0.45) {$x$-direction};
        \node[font=\scriptsize, black!70, rotate=90] at (-0.55,2.8) {$y$-direction};

        \draw[-{Stealth[length=2.0mm]}, black!55] (2.2,-0.35) -- (3.55,-0.35);
        \draw[-{Stealth[length=2.0mm]}, black!55] (-0.35,1.95) -- (-0.35,3.1);

        \node[font=\scriptsize, black!75] at (0,4.55)
            {$\text{circular CoP} \;=\; \{\,p(\theta): 0\le \theta < 2\pi\,\}$};
    \end{scope}
\end{tikzpicture}
\caption{Bicycle-wheel geometry of the CoP: a circular body with spokes as axes and rim points as the terminal locations of points in the structure.}
\end{figure}

\begin{example}
The following are some examples of CoPs with prime number base sets $\mathbb{P}$
$$
\mathcal{C}(36,\mathbb{P})=\lbrace [5],\mathbf{[7]},[13],[17],\mathbf{[19]},[23],[29],\mathbf{[31]}\rbrace
$$
and 
$$
\mathcal{C}(38,\mathbb{P})=\lbrace [7],[19],[31]\rbrace.
$$ 
with their corresponding weight sets 
$$
||\mathcal{C}(36,\mathbb{P})||:=\{5,7,13,17,19,23,29,31\}
$$ 
and 
$$
||\mathcal{C}(38,\mathbb{P})||:=\{7,19,31\}.
$$
\end{example}
\bigskip

The binary Goldbach and Lemoine conjecture can be reformulated in the language of circle of partitions as

\begin{conjecture}[The binary Goldbach conjecture]
Let $\mathbb{P}$ and $2\mathbb{N}$ denote the set of all prime numbers and the set of all positive even numbers, respectively. For all $n\in 2\mathbb{N}$ with $n\geq 6$, we have
$$
\mathcal{C}(n,\mathbb{P})\neq \emptyset.
$$
\end{conjecture}

\begin{conjecture}[Lemoine conjecture]
Let $\mathbb{P}$ and $2\mathbb{P}$ denote the set of all prime numbers and their doubles, respectively. For all $n\in 2\mathbb{N}+1$ with $n\geq 7$, we have
$$
\mathcal{C}(n,\mathbb{P}\cup 2\mathbb{P})\neq \emptyset.
$$
\end{conjecture}
\bigskip

It is clear that the two preceding statements are the precise reformulation of the binary Goldbach and Lemoine conjectures in the language of the circle of partitions. Indeed in any of these cases, showing that the circle of partition is always non-empty for each generator within the specified range would mean that there must exist a point in the circle of partition. Consequently, it will follow from the properties of circle of partitions that each point must have an axes partner. In other words, if a circle of partition is non-empty, then it must have at least an axes.
\bigskip

\section{The Squeeze Principle}\label{sec:squeeze principle}

In this section, we introduce the \emph{squeeze principle} and its consequences if the set of all odd prime numbers is the base set of CoPs.

\begin{theorem}[The squeeze principle]\label{L_squeeze principle}
Let $\B\subset\M\subseteq\N$, and suppose that $\C(m,\B)$ and $\C(m+t,\B)\neq\emptyset$ for some $t\geq 4$. Suppose that there exist $\mathbb{L}_{[x],[y]}\inn \mathcal{C}(m+t,\mathbb{M})$ with $x\in \mathbb{B}$ and $x<y$ such that 
\begin{align}\label{E_maxcovered}
y>w= \max\set{u\in||\C(m,\M)||\mid u\in \B}>m-x.
\end{align}
There exists a generator $s$ with $m<s<m+t$ such that $\C(s,\B)\neq\emptyset$.
\end{theorem}

\begin{proof}
By \eqref{E_maxcovered}, we have $w\in\B$. By hypothesis, the axis $\mathbb{L}_{[x],[y]}\inn \mathcal{C}(m+t,\mathbb{M})$ exists with $x\in \mathbb{B}$ such that $m-w<x<y$. We deduce
\begin{align}
m=w+(m-w)<\ul{w+x}&=w+(m+t-y)=m+t+(w-y)\nonumber\\
&<m+t \mbox{, since }y>w \label{E_inequality}
\end{align}
and $m-w<x=m+t-y$ implies $y-w<t$. With $s=\ul{w+x}$, there is an axis $\LL_{[x],[w]}\inn\C(s,\B)$ and it follows that $\mathcal{C}(s,\mathbb{B})\neq \emptyset$ with $m<s<m+t$. 
\end{proof}
\bigskip

Theorem \ref{L_squeeze principle} may be viewed as a basic apparatus for studying the possibility of partitioning numbers of a particular parity into components that belong to a special subset of positive integers. It works by choosing two non-empty CoPs with the same base set and finding further non-empty CoPs with generators trapped in between these two generators. This principle can be used with care to study the broader question that concerns the feasibility of partitioning numbers with each summand belonging to the same subset of positive integers. We state the following proposition as an outgrowth of Theorem \ref{L_squeeze principle}.

\begin{proposition}[The interval binary Goldbach partition detector]\label{P_The interval binary Goldbach partition detector}
Let $\mathbb{P}$ be the set of all prime numbers and let $\mathcal{C}(m,\mathbb{P}),\mathcal{C}(m+t,\mathbb{P})\neq \emptyset$ for some $t\geq 4$. Suppose that there exist $\mathbb{L}_{[x],[y]}\inn \mathcal{C}(m+t,\mathbb{N})$ with $x\in\p$ and $x<y$ such that
\begin{align}\label{E_maxcoveredP}
y>w= \max\set{u\in||\C(m,\N)||\mid u\in \p}>m-x.
\end{align}
There exists a generator $s$ with $m<s<m+t$ such that $\C(s,\p)\neq \emptyset$.
\end{proposition}

\begin{proof}
This is a simple consequence of Theorem \ref{L_squeeze principle} by taking $\M=\N$ and $\B=\p$.
\end{proof}

\begin{proposition}[Interval Goldbach partition]\label{P_intervalGoldbach}
Let $\p$ be the set of all prime numbers and $\C(m,\p), \C(m+t,\p)\neq\emptyset$ be for some $t\geq 4$. If $m-1\in \p$, then there exists some $s\equiv 0\pmod 2$ with $m<s<m+t$ such that $\C(s,\p)\neq\emptyset$.
\end{proposition}

\begin{proof}
By the requirements $\C(m,\p),\C(m + t,\p)\neq\emptyset$ for $t \geq 4$ and with $w$ by virtue of \eqref{E_maxcoveredP}, we choose
$\LL_{[3],[y]}\inn \C(m + t,\N)$ so that $w=m-1$ and $y>w$
since $y=m+t-3>m$ for $t\geq 4$ and $m-1\in\p$. We deduce
\begin{align*}
y-w=y-(m-1)\leq(m+t-3)-(m-1)<t
\end{align*}
and the conditions in Proposition \ref{P_The interval binary Goldbach partition detector} are satisfied, so that there exist some $s\equiv 0\pmod 2$ with $m<s<m+t$ such that $\C(s, \p)\neq\emptyset$, f.i. $s=3+m-1=m+2$ with $\LL_{[3],[m-1]}\inn\C(m+2,\p)$.
\end{proof}

\begin{proposition}\label{P_finite Goldbach in an interval}
Let $\mathbb{P}$ be the set of all prime numbers and $\mathcal{C}(m,\mathbb{P}),\mathcal{C}(m+t,\mathbb{P})\neq \emptyset$ be for some $t\geq 4$ such that $m-1\in \mathbb{P}$. There are finitely many $s\equiv 0\pmod 2$ with $m<s<m+t$ such that $\C(s,\p)\neq \emptyset$.
\end{proposition}

\begin{proof}
The result is obtained by iterating on the generators $s\equiv 0\pmod 2$ with $m<s<m+t$ such that $\C(s,\p)\neq \emptyset$.
\end{proof}

\begin{theorem}[Conditional Goldbach]\label{T_Conditional Goldbach}
Let $\mathbb{P}$ be the set of all prime numbers and $m\in 2\mathbb{N}$ such that $\mathcal{C}(m,\mathbb{P})\neq \emptyset$ for a \textbf{sufficiently} large $m$ . If for \textbf{all} $t\geq 4$ there exists $\mathbb{L}_{[x],[y]}\inn \mathcal{C}(m+t,\mathbb{N})$ with $x\in \mathbb{P}$ and $x<y$ such that
\[
y> w= \max\set{u\in||\C(m,\N)||\mid u\in \p}>m-x,
\]
then there are CoPs $\mathcal{C}(s,\mathbb{P})\neq \emptyset$ for all (sufficiently large) $s\in 2\mathbb{N}$ with $s>m$.
\end{theorem}

\begin{proof}
It is known that there are infinitely many even numbers that can be written as the sum of two primes, so that for \textbf{sufficiently} large $m\in 2\mathbb{N}$ with $\mathcal{C}(m,\mathbb{P})\neq \emptyset$, then $t\geq 4$ can be \textbf{arbitrarily} chosen large such that $\mathcal{C}(m+t,\p) \neq \emptyset$. By the requirements and Proposition \ref{P_The interval binary Goldbach partition detector}, there must exist some $s\equiv 0\pmod 2$ with $m<s<m+t$ such that $\mathcal{C}(s,\mathbb{P})\neq \emptyset$. Now, we continue our arguments on the intervals of the generators $[m,s]$ and $[s,s+r]$. If there exist some $u,v\in 2\mathbb{N}$ such that $m<u<s$ and $s<v<s+r$, then we repeat the argument under the requirements (for arbitrary $t$) to deduce that $\mathcal{C}(u,\mathbb{P})\neq \emptyset$ and $\mathcal{C}(v,\mathbb{P})\neq \emptyset$. We can iterate the process so long as there exist some even generators trapped in the following sub-intervals of generators $[m,u],[u,s],[s,v],[v,v+r]$ where $v+r=m+t$ for $t\geq 4$. Since $t$ can be arbitrarily chosen so that $\mathcal{C}(m+t,\p) \neq \emptyset$, the assertion follows immediately.
\end{proof}
\bigskip

Now, we use the \emph{squeeze principle} to solve the Lemoine conjecture in analogy to its use for the binary Goldbach conjecture above.

\begin{proposition}[The first interval Lemoine partition detector]\label{P_The first interval Lemoine partition detector}
Let $\mathbb{P}$ and $2\mathbb{P}$ be the set of all prime numbers and their doubles, respectively, and let $\mathcal{C}(m,\mathbb{P}\cup 2\mathbb{P}),\mathcal{C}(m+t,\mathbb{P}\cup 2\mathbb{P})\neq \emptyset$ be for some $t\geq 4$. Suppose that there exist $\mathbb{L}_{[x],[y]}\inn \mathcal{C}(m+t,\mathbb{N})$ with $x\in\p$ and $x<y$ such that
\begin{align}\label{E_maxcoveredP2}
y>w=\max\set{u\in||\C(m,\N)||\mid u\in \p\cup 2\p}\in 2\mathbb{P}>m-x.
\end{align}
There exists a generator $s$ with $m<s<m+t$ such that $\C(s,\p\cup 2\p)\neq \emptyset$.
\end{proposition}

\begin{proof}
This is a consequence of Theorem \ref{L_squeeze principle} when we take $\M=\N$ and $\B=\p\cup 2\p$.
\end{proof}

\begin{proposition}[The second interval Lemoine partition detector]\label{P_The second interval Lemoine partition detector}
Let $\mathbb{P}$ and $2\mathbb{P}$ be the set of all prime numbers and their doubles, respectively, and $\mathcal{C}(m,\mathbb{P}\cup 2\mathbb{P}),\mathcal{C}(m+t,\mathbb{P}\cup 2\mathbb{P}) \neq \emptyset$ by $t\geq 4$. Suppose that there exist $\mathbb{L}_{[x],[y]}\inn \mathcal{C}(m+t,\mathbb{N})$ with $x\in 2\p$ and $x<y$ such that
\begin{align}\label{E_maxcoveredP3}
y>w= \max\set{u\in||\C(m,\N)||\mid u\in \p\cup 2\p}\in \mathbb{P}>m-x.
\end{align}
There exists a generator $s$ with $m<s<m+t$ such that $\C(s,\p\cup 2\p)\neq \emptyset$.
\end{proposition}

\begin{proof}
The proof is the same as in Proposition \ref{P_The first interval Lemoine partition detector}.
\end{proof}

\begin{theorem}[Conditional Lemoine]\label{T_Conditional Lemoine}
Let $\mathbb{P}$ and $2\mathbb{P}$ be the set of all prime numbers and their doubles, respectively, and $m\in 2\mathbb{N}+1$ be such that $\mathcal{C}(m,\mathbb{P})\neq \emptyset$ for a \textbf{sufficiently} large $m$. If for \textbf{all} $t\geq 4$ there exists $\mathbb{L}_{[x],[y]}\inn \mathcal{C}(m+t,\mathbb{N})$ with $x\in \mathbb{P}$ and $x<y$ such that
\[
y>w=\max\set{u\in||\C(m,\N)||\mid u\in \p\cup 2\p}\in 2\p>m-x,
\]
or there exist $\mathbb{L}_{[x],[y]}\inn \mathcal{C}(m+t,\mathbb{N})$ with $x\in 2\mathbb{P}$ and $x<y$ such that
\[
y> w= \max\set{u\in||\C(m,\N)||\mid u\in \p\cup 2\p}\in \p>m-x
\]
then there are CoPs $\mathcal{C}(s,\mathbb{P}\cup 2\mathbb{P})\neq \emptyset$ for all (sufficiently large) $s\in 2\mathbb{N}+1$ with $s>m$.
\end{theorem}

\begin{proof}
It is known that there are infinitely many odd numbers that can be written as the sum of a prime and a double of a prime, so that for a \textbf{sufficiently} large $m\in 2\mathbb{N}+1$ with $\mathcal{C}(m,\mathbb{P}\cup 2\mathbb{P})\neq \emptyset$, then $t\geq 4$ can be \textbf{arbitrarily} chosen large such that $\mathcal{C}(m+t,\p)\neq \emptyset$. By the requirements and Propositions \ref{P_The first interval Lemoine partition detector} and \ref{P_The second interval Lemoine partition detector}, there must exist some $s\equiv 1\pmod 2$ with $m<s<m+t$ such that $\mathcal{C}(s,\mathbb{P}\cup 2\mathbb{P})\neq \emptyset$. Now, we continue our arguments on the intervals of the generators $[m,s]$ and $[s,s+r]$. If there exist some $u,v\in 2\mathbb{N}+1$ such that $m<u<s$ and $s<v<s+r$, then we repeat the argument under the requirements (for arbitrary $t$) to deduce that $\mathcal{C}(u,\mathbb{P}\cup 2\mathbb{P})\neq \emptyset$ and $\mathcal{C}(v,\mathbb{P}\cup 2\mathbb{P})\neq \emptyset$. We can iterate the process so long as there exist some odd generators (numbers) trapped in the following sub-intervals of generators $[m,u],[u,s],[s,v],[v,v+r]$ where $v+r=m+t$ for $t\geq 4$. Since $t$ can be arbitrarily chosen so that $\mathcal{C}(m+t,\p\cup 2\p)\neq \emptyset$, the assertion follows immediately.
\end{proof}
\bigskip

\section{Application to the Binary Goldbach Conjecture}\label{sec:application to Goldbach}

In this section, we apply the \emph{squeeze principle} to study the asymptotic version of the binary Goldbach conjecture. Despite the Estermann proof from 1938 (see, e.g, \cite{estermann1938goldbach}, \cite{chudakov1938goldbach}) that the binary Goldbach conjecture is true for almost all positive integers, and the work of Chen that every sufficiently large even number can be written as a prime and almost prime (see, e.g, \cite{chen2002representation}) we can use our elementary tool to independently establish and prove the binary Goldbach conjecture in an asymptotic sense. We gather the following established facts that will feature prominently in our arguments.

\begin{lemma}[The prime number theorem]\label{L_the prime number theorem}
Let $\pi(m)$ denote the number of prime numbers less than or equal to $m$ and $p_{\pi(m)}$ denote the $\pi(m)^{th}$ prime number. We have the asymptotic relation
\begin{align}
p_{\pi(m)}\sim \pi(m)\log \pi(m) \quad \text{and} \quad \pi(m)\sim \frac{m}{\log m}.\nonumber
\end{align}
\end{lemma}

In keeping with the notation of the previous section, we write
\begin{align}\label{E_w-max-p}
w=\max\set{u\in||\C(m,\N)|| \mid u\in\p}=p_{\pi(m)}.
\end{align}
and set
$$
w''=2p_{\pi(\frac{m}{2})}.
$$

\bigskip

\begin{lemma}[Bertrand's postulate]\label{L_Bertrand's postulate}
There exists a prime number in the interval $(k,2k)$ for all $k>1$.
\end{lemma}
\bigskip

\begin{lemma}[The little lemma]\label{L_the little lemma}
Let $\mathbb{P}$ be the set of all prime numbers and $m\in \mathbb{N}$ be \textbf{sufficiently} large such that $\mathcal{C}(m,\mathbb{P})\neq \emptyset$. There exists some $h:=h(m)=o(\log m)$ such that
\begin{align*}
m-w''\sim m-w\sim h(m)\frac{m}{\log m}.
\end{align*}
For all $x\in \mathbb{P}$ satisfying 
$$
m-w<x\leq 3(m-w)
$$ 
the inequality
\begin{align*}
0\lesssim |w-(m+t-x)|\lesssim t
\end{align*}
hold for $t\geq 4$.
\end{lemma}

\begin{proof}
The first assertion can be deduced from the prime number theorem. By the Bertrand postulate, there exists an odd prime number $x\in \mathbb{P}$ satisfying $m-w<x\leq 3(m-w)$, since 
\begin{align}
w=\max\set{u\in||\C(m,\N)|| \mid u\in\p}=p_{\pi(m)}.
\end{align}
imply $m-w\geq 1$.
We deduce
\begin{align*}
m+t-x&\gtrsim m+t-3h(m)\frac{m}{\log m}\\
&=m\left(1-3\frac{h(m)}{\log m}\right)+t\\
&\sim m+t>p_{\pi(m)}=w
\end{align*}
and 
\begin{align*}
|w-(m+t-x)|&=|m+t-x-w|\\
&<|m+t-(m-w)-w|\\
&=t\nonumber
\end{align*}
for $t\geq 4$.
\end{proof}
\bigskip

We are now ready to prove the binary Goldbach conjecture for all \textbf{sufficiently} large even numbers. The following result is a combination of ideas developed in this paper. 

\begin{theorem}[Asymptotic Goldbach theorem]\label{T_Asymoptotic Goldbach theorem}
Every \textbf{sufficiently} large even number can be written as the sum of two prime numbers.
\end{theorem}

\begin{proof}
The claim is equivalent to the following statement:
\begin{center}
For every sufficiently large even number $n$, we have $\C(n,\p)\neq\emptyset$.
\end{center}
It is known that there are infinitely many even numbers $m>0$ with $\mathcal{C}(m,\mathbb{P})\neq \emptyset$. We choose a \textbf{sufficiently} large $m\in 2\N$ such that $\mathcal{C}(m,\mathbb{P})\neq \emptyset$ and choose $t\geq 4$ such that $\mathcal{C}(m+t,\mathbb{P})\neq \emptyset$. We set 
\begin{align}
w=\max\set{u\in||\C(m,\N)|| \mid u\in\p}=p_{\pi(m)}\nonumber
\end{align}
and choose a prime number $x\leq 3(m-w)$ such that $x>m-w$, since by the Bertrand postulate (Lemma \ref{L_Bertrand's postulate}) there exists a prime number $x$ such that $x\in (k,2k)$ for every $k>1$. We obtain for the axis partner $[y]$ of the axis point $[x]$ of $\LL_{[x],[y]}\inn\C(m+t,\N)$ the inequality 
\begin{align*}
y=m+t-x&\gtrsim m+t-3h(m)\frac{m}{\log m}\\
&=m\left(1-3\frac{h(m)}{\log m}\right)+t\\
&\sim m+t>p_{\pi(m)}=w
\end{align*}
for $t\geq 4$ and by the lemma \ref{L_the little lemma}, we deduce
\begin{align*}
h(m)\frac{m}{\log m}\sim m-w<x
\end{align*}
with $h(m)=o(\log m)$ and 
\begin{align*}
\abs{y-w}=\abs{(m+t-x)-w}=\abs{m-w+t-x}\lesssim\abs{x+t-x}= t.
\end{align*}
The requirements in Theorem \ref{T_Conditional Goldbach} are \textbf{asymptotically} fulfilled 
with
\[
y\gtrsim w\mbox{ and }x\gtrsim m-w\mbox{ and }0\lesssim \abs{y-w}\lesssim t.
\]
The result follows by arbitrarily choosing $t\geq 4$ so that $\mathcal{C}(m+t,\mathbb{P})\neq \emptyset$ and adapting the proof in Theorem \ref{T_Conditional Goldbach}.
\end{proof}
\bigskip

\section{Application to the Lemoine Conjecture}\label{sec:application to Lemoine}

In this section, we apply the \emph{squeeze principle} to study Lemoine conjecture (see, e.g, \cite{levi1963}, \cite{lemoine1894}) in the asymptotic.

\begin{lemma}\label{L_the prime number theorem2}
Let $\pi(m)$ denote the number of prime numbers less than or equal to $m$ and $p_{\pi(m)}$ denote the $\pi(m)^{th}$ prime number. We have the asymptotic relation
\begin{align*}
p_{\pi(m)}\sim 2p_{\pi(\frac{m}{2})}.
\end{align*}
\end{lemma}

\begin{proof}
This is the prime number theorem.
\end{proof}

In keeping with the notation of the previous section, we write
\begin{align}\label{E_w-max-p1}
w=\max\set{u\in||\C(m,\N)||~\mid u\in\p\cup 2\p}=p_{\pi(m)}
\end{align}
provided that $w\in \mathbb{P}$ and 
\begin{align}\label{E_w-max-p2}
w'=\max\set{u\in||\C(m,\N)||~\mid u\in\p\cup 2\p}=2p_{\pi(\frac{m}{2})}
\end{align}
provided that $w'\in 2\mathbb{P}$.
\bigskip

\begin{lemma}[The first little lemma]\label{L_the first little lemma}
Let $\mathbb{P}$ and $2\mathbb{P}$ be the set of all prime numbers and their doubles, respectively, and let $m\in \mathbb{N}$ be \textbf{sufficiently} large such that $\mathcal{C}(m,\mathbb{P}\cup 2\mathbb{P})\neq \emptyset$. Suppose that
\begin{align}
w=\max\set{u\in||\C(m,\N)||~\mid u\in\p\cup 2\p}=p_{\pi(m)}.
\end{align}
There exists some $h:=h(m)=o(\log m)$ such that the asymptotic relation
\begin{align*}
m-w\sim m-w'\sim h(m)\frac{m}{\log m}
\end{align*}
holds. Furthermore, for all $x'\in 2\mathbb{P}$ with $x'=2x$ for $x\in \mathbb{P}$ satisfying 
\[
\frac{1}{2}(m-w')<x<m-w'
\]
the inequality
\begin{align*}
0\lesssim |w-(m+t-x')|\lesssim t
\end{align*}
holds for $t\geq 4$.
\end{lemma}

\begin{proof}
The first assertion can be deduced from the prime number theorem. By the lemma \ref{L_Bertrand's postulate}, there is a prime between $\frac{1}{2}(m-w')$ and $m-w'$, since 
\begin{align}
w=\max\set{u\in||\C(m,\N)||~\mid u\in\p\cup 2\p}=p_{\pi(m)}.
\end{align}
imply $m-w'\geq 3$.
We deduce
\begin{align*}
m+t-x'=m+t-2x&\gtrsim m+t-2h(m)\frac{m}{\log m}\\
&=m\left(1-2\frac{h(m)}{\log m}\right)+t\\
&\sim m+t\geq p_{\pi(m)}=w
\end{align*}
and 
\begin{align*}
|w-(m+t-x')|&\sim |m+t-x'-w'|\\
&\lesssim |m+t-(m-w')-w'|\\
\\
&=t\nonumber
\end{align*}
for $t\geq 4$.
\end{proof}

\begin{lemma}[The second little lemma]\label{L_the second little lemma}
Let $\mathbb{P}$ be the set of all prime numbers and their doubles, respectively, and $m\in \mathbb{N}$ be \textbf{sufficiently} large such that $\mathcal{C}(m,\mathbb{P}\cup 2\mathbb{P})\neq \emptyset$. Suppose that 
\begin{align}
w'=\max\set{u\in||\C(m,\N)||~\mid u\in\p\cup 2\p}=2p_{\pi(\frac{m}{2})}.
\end{align}
There exists some $h:=h(m)=o(\log m)$ such that 
$$
m-w'\sim m-w\sim h(m)\frac{m}{\log m}.
$$
For all $x\in \mathbb{P}$ satisfying 
\[
m-w<x<2(m-w)
\]
we have
\begin{align*}
0\lesssim |w'-(m+t-x)|\lesssim t
\end{align*}
for $t\geq 4$.
\end{lemma}

\begin{proof}
By the lemma \ref{L_Bertrand's postulate}, there is a prime between $m-w$ and $2(m-w)$, since 
\begin{align}
w'=\max\set{u\in||\C(m,\N)||~\mid u\in\p\cup 2\p}=2p_{\pi(\frac{m}{2})}
\end{align}
imply $m-w\geq 2$. We deduce
\begin{align*}
m+t-x&>m+t-2h(m)\frac{m}{\log m}\\
&=m\left(1-2\frac{h(m)}{\log m}\right)+t\\
&\sim m+t\geq p_{\pi(m)}\sim 2p_{\pi(\frac{m}{2})}=w'
\end{align*}
and 
\begin{align*}
|w'-(m+t-x)|&=|m+t-x-w'|\\
&\lesssim |m+t-(m-w')-w'|\\
&=t\nonumber
\end{align*}
for $t\geq 4$.
\end{proof}
\bigskip

We are now ready to prove the Lemoine conjecture for all \textbf{sufficiently} large odd numbers. It is a case-by-case argument and a culmination of ideas developed in this paper. 

\begin{theorem}[Asymptotic Lemoine theorem]\label{T_Asymoptotic Lemoine theorem}
Every \textbf{sufficiently} large odd number can be written as the sum of a prime number and a double of a prime number.
\end{theorem}

\begin{proof}
The claim is equivalent to the following statement:
\begin{center}
For every sufficiently large odd number $n\in 2\mathbb{N}+1$, we have $\C(n,\p\cup 2\p)\neq\emptyset$, since only the sum of an odd and an even number provides an odd number and therefore each axis $\LL_{[x],[y]}\inn\C(m,\p\cup 2\p)$ has an odd and an even axis point.
\end{center}
It is known that there are infinitely many odd numbers $m>0$ with $\mathcal{C}(m,\mathbb{P}\cup 2\mathbb{P})\neq \emptyset$. We choose a \textbf{sufficiently} large $m\in 2\N+1$ such that $\mathcal{C}(m,\mathbb{P}\cup 2\mathbb{P})\neq \emptyset$ and choose $t\geq 4$ such that $\mathcal{C}(m+t,\mathbb{P}\cup 2\mathbb{P})\neq \emptyset$. Now, we distinguish and examine two special cases as follows:
\begin{enumerate}
\item  The case $$w=\max\set{u\in||\C(m,\N)||~\mid u\in\p\cup 2\p}=p_{\pi(m)}$$

\item The case $$w'=\max\set{u\in||\C(m,\N)||~\mid u\in\p\cup 2\p}=2p_{\pi(\frac{m}{2})}$$
\end{enumerate} 
In the case 
$$
w=\max\set{u\in||\C(m,\N)||~\mid u\in\p\cup 2\p}=p_{\pi(m)}
$$ 
we choose a prime number $x<m-w'$ such that $x>\frac{1}{2}(m-w')$, since by the Bertrand postulate (Lemma \ref{L_Bertrand's postulate}) there exists a prime number $x$ such that $x\in (k,2k)$ for every $k>1$ and set $2x=x'\in 2\mathbb{P}$. We obtain for the axis partner $[y']$ of the axis point $[x']$ of $\LL_{[x'],[y']}\inn\C(m+t,\N)$ the inequality 
\begin{align*}
y'=m+t-x'&\gtrsim m+t-2h(m)\frac{m}{\log m}\\
&=m\left(1-2\frac{h(m)}{\log m}\right)+t\\
&\sim m+t\geq p_{\pi(m)}=w
\end{align*}
for $t\geq 4$ and by the lemma \ref{L_the first little lemma} the following asymptotic inequalities
\begin{align*}
m-w\sim m-w'<x'
\end{align*}
and 
\begin{align*}
\abs{y'-w}=\abs{(m+t-x')-w}=\abs{m-w+t-x'}\lesssim\abs{x'+t-x'}= t.
\end{align*}
The requirements in Theorem \ref{T_Conditional Lemoine} are \textbf{asymptotically} fulfilled in this case
with
\[
y'\gtrsim w\mbox{ and }x'\gtrsim m-w\mbox{ and }0\lesssim \abs{y'-w}\lesssim t.
\]
In the case 
$$
w'=\max\set{u\in||\C(m,\N)||~\mid u\in\p\cup 2\p}=2p_{\pi(\frac{m}{2})}
$$ 
we choose a prime number $x<2(m-w)$ such that $x>(m-w)$, since by Bertrand's postulate (Lemma \ref{L_Bertrand's postulate}) there exists a prime number $x$ such that $x\in (k,2k)$ for every $k>1$. We obtain for the axis partner $[y]$ of the axis point $[x]$ of $\LL_{[x],[y]}\inn\C(m+t,\N)$ the inequality 
\begin{align*}
y=m+t-x&\gtrsim m+t-2h(m)\frac{m}{\log m}\\
&=m\left(1-2\frac{h(m)}{\log m}\right)+t\\
&\sim m+t\geq p_{\pi(m)}\sim 2p_{\pi(\frac{m}{2})}=w'
\end{align*}
for $t\geq 4$ and by appealing to Lemma \ref{L_the second little lemma} the following asymptotic inequalities
\begin{align*}
h(m)\frac{m}{\log m}\sim m-w'\sim m-w<x
\end{align*}
where $h(m):=o(\log m)$ and 
\begin{align*}
\abs{y-w'}=\abs{(m+t-x)-w'}=\abs{m-w'+t-x}\lesssim\abs{x+t-x}= t.
\end{align*}
The requirements in Theorem \ref{T_Conditional Lemoine} are \textbf{asymptotically} fulfilled in this second case
with
\[
y\gtrsim w'\mbox{ and }x\gtrsim m-w'\mbox{ and }0\lesssim \abs{y-w'}\lesssim t.
\] 
The result follows by arbitrarily choosing $t\geq 4$ so that $\mathcal{C}(m+t,\mathbb{P}\cup 2\mathbb{P})\neq \emptyset$ and adapting the proof in Theorem \ref{T_Conditional Lemoine}.
\end{proof}
\section{Further remarks}\label{sec:conclusion}

Theorem \ref{T_Asymoptotic Goldbach theorem} and Theorem \ref{T_Asymoptotic Lemoine theorem} are both equivalent to the statement: there must exist some $N>0$ such that for all $m\geq N$ every even number $m$ can be written as the sum of two prime numbers (resp. every odd number is the sum of a prime and a double of a prime). This result--although constructive to some extent--loses its constructive flavour in a way that one cannot perform this construction to cover all odd numbers, since we are unable to obtain any quantitative (lower) bound for the threshold $N$. Regardless, we are able to \textbf{asymptotically} get a handle on the conjecture. The asymptotic version of the Lemoine conjecture implies an asymptotic version ternary of the Goldbach conjecture \cite{helfgott2013ternary}.


\end{document}